\documentclass[12pt]{article}
\usepackage{amsmath, amsthm, amssymb}
\usepackage{amsfonts}
\usepackage{epsfig}
\usepackage{graphics}
\usepackage{subfigure}
\usepackage{graphicx}
\usepackage{color}
\usepackage{epstopdf} 
\usepackage{float} 
\usepackage{mathrsfs} 
\usepackage{relsize} 

\counterwithin*{equation}{section}
\newcommand{\const}{{\rm const}}
\newcommand{\ds}{\displaystyle}

\newtheorem{thm}{Theorem}
\newtheorem{df}{Definition}
\newtheorem{lm}{Lemma}
\newtheorem{rem}{Remark}
\newcommand{\bq}{\begin{equation}}
\newcommand{\eq}{\end{equation}}
\newcommand{\bqr}{\begin{eqnarray}}
\newcommand{\eqr}{\end{eqnarray}}
\newcommand{\bqrn}{\begin{eqnarray*}}
\newcommand{\eqrn}{\end{eqnarray*}}

\title{Fractional Quasi-Bessel Equations} 
\author{\normalsize Pavel B. Dubovski $^1$, Jeffrey A. Slepoi $^1$}

\begin{document}
\begin{abstract}
	
In this paper we consider fractional quasi-Bessel equations $$\sum_{i=1}^{m}d_i x^{\alpha_i+p_i}D^{\alpha_i} u(x) + (x^\beta - \nu^2)u(x)=0$$ and construct their existence and uniqueness theory in the class of fractional series. Our methodology allows us to obtain new results for a broad class of fractional differential equations including Cauchy-Euler and constant-coefficient equations. 
	  \medskip
	 
	 
	 \smallskip
	 
	 
 \end{abstract}
\date{ }
\maketitle

\section{Introduction}
We generalize the classical Bessel equation \eqref{BesselEq} 
\begin{equation}\label{BesselEq}
	x^2 u'' + x u' + (x^2 - \nu^2)u =0
\end{equation}
to a broad class of equations including multi-term fractional Bessel, Cauchy-Euler equations, and fractional differential equations with constant coefficients.
\noindent We consider the Caputo  
\begin{equation}\label{CaputoFrD}
	D_C^{\alpha}u(t)=\frac{1}{\Gamma(n-\alpha)}\int_{0}^{t}\frac{u^{(n)}(x) dx}{(t-x)^{\alpha+1-n}},n-1 \le \alpha < n, n \in \mathbb{N} 
\end{equation}
and Riemann-Liouville 
\begin{equation}\label{RLFrD}
	D_{R}^{\alpha}u(t)=\frac{1}{\Gamma(n-\alpha)}\frac{d^{n}}{dt^n}\int_{0}^{t}\frac{u(x) dx}{(t-x)^{\alpha+1-n}}, n-1 \le \alpha < n, n \in \mathbb{N}
\end{equation}
fractional derivatives.\\

Equation (\ref{BesselEq}) is a source for many generalizations and extensions. Numerous results for  hyper-Bessel equations and corresponding operators, introduced by Dimovski, are obtained in \cite{2018Al-Musalhi, 
 2021Droghei, Kir1994, 2008Kiryakova}. 
Paper \cite{2015Bengochea} introduces another definition for the fractional order Bessel operator, which allows to construct an operational calculus applicable to Bessel equation with fractional derivatives.
W. Okrasiński and L. Plociniczak \cite{2013Okrasinski} consider a fractional modification of the Bessel equation, and solve it in terms of the power series.

A natural extension of the classical Bessel equation (\ref{BesselEq}) in terms of Caputo fractional derivatives 
\begin{equation}\label{BesselFrDE}
	x^{2\alpha}D^{2\alpha}u(x)+x^\alpha D^\alpha u(x)+(x^{2\alpha}-\nu^2)u(x)=0, \alpha \in (0,1]
\end{equation}
was analyzed by Rodrigues, Viera and Yakubovich \cite{Rodrigues}, where a solution in a form of series for equation \eqref{BesselFrDE} was identified for a some specific values of $\nu$ depending on $\alpha$.\\

The multi-term fractional Bessel equation  
\begin{equation}\label{GenBesselEqn}
\sum_{i=1}^{m_1}d_i x^{\alpha_i}D_C^{\alpha_i} u(x) + (x^\beta - \nu^2)u(x)=0, \alpha_i > 0, \beta > 0
\end{equation}
with existence and uniqueness criteria was presented by the authors in \cite{2020DuSl-3}
where we constructed solutions in the form of series
\begin{equation}\label{SolnFn}
	u(x)=\sum_{n=0}^{\infty}c_n x^{\gamma+\beta n}.  
\end{equation}
with coefficients $c_n$ as 
\begin{eqnarray} \label{Eqnforc_n}
	c_{n}
	=\frac{(-1)^n}{\ds\prod_{k=1}^{n}\left(\ds\sum_{i=1}^{m_1}\frac{d_i\Gamma(1+\gamma+\beta k)}{\Gamma(1+\gamma+\beta k-\alpha_i)}-\nu^2\right)}.
\end{eqnarray}\\

The characteristic equation for equation \eqref{GenBesselEqn}
\begin{equation}\label{CondForC0}
	\sum_{i=1}^{m_1}\frac{d_i\Gamma(1+\gamma)}{\Gamma(1+\gamma-\alpha_i)}-\nu^2=0
\end{equation}
must be solvable and all derivatives $D^{\alpha_i}x^\gamma$ should exist. 
%
%
Theorem \ref{Thrm1} below states that series solution (\ref{SolnFn}),(\ref{Eqnforc_n}) always exists in real numbers as long as all $d_i > 0$ and $\nu$ satisfies a certain threshold condition.
Our results for equation \eqref{GenBesselEqn} are obtained for a broad class of parameters $\alpha_i$ and $\nu$. They include some results of (\cite{Kir1994}, sec. 3.4) and \cite{Rodrigues}. For example, if we deal with derivatives of integer orders $\alpha_i$ and $\nu=0$ in equation (\ref{GenBesselEqn}), then we arrive at the eigenvalue problem, 
which is solved in \cite{Kir1994} and where the fundamental system of solutions (in terms of the hyper-Bessel functions of Delerue) has been determined. Besides fractional Bessel equations (\ref{GenBesselEqn}), the quasi-Bessel equations generalize 
Cauchy-Euler equidimensional equations as well and, thus, we extend further a number of results from the works of Kilbas and Zhukovskaya \cite{2009Kilbas-Zh, 2011Zh-Kilbas, 2018Zhukovskaya}. After a minor modification, our approach becomes also 
applicable to the equations with constant coefficients. In works \cite{2014Kim}, using the methods of \cite{2009Hilfer}, the initial value problem 
is solved for constant-coefficient fractional equations. In \cite{2014Atanackovic}, the Cauchy problems for some linear constant-coefficient equations with real, including irrational, derivatives, are analyzed and solved. The presented method and its particular version for constant-coefficient equations also allow to treat a class of irrational derivatives.
Additional extensive information regarding fractional differential equations can be found in many works on fractional calculus, e.g., in monographs \cite{Kilbas, Kir1994, Podlubny, Samko}.\\

In this paper we analyze the next generalization of multi-term Bessel equations -- the quasi-Bessel fractional equations
\begin{equation}\label{QuasiBessel}
\sum_{i=1}^{m}d_i x^{\alpha_i+p_i}D^{\alpha_i} u(x) + (x^\beta - \nu^2)u(x)=0, \ \ x > 0,
\end{equation}
where $\alpha_i \in \mathbb{R^+}$ with $\alpha_1 = \max\limits_{1\leq i\leq m}\{\alpha_i\}$ and $p_1=0$. Shifting indices $p_i$ (deviations with respect to (\ref{GenBesselEqn}))
should be non-negative and satisfy special almost-rational conditions defined in Section \ref{SectBeyondFracBesEq}. 
If all additional indices $p_i$ are zero, then we arrive at the previously investigated in \cite{2020DuSl-3} multi-term fractional Bessel equation (\ref{GenBesselEqn}), where the powers of multipliers $x$ must match the order of derivatives. However, in this research the powers of $x$ have to match the order of only one derivative -- the derivative of the highest order. As we show below, many other equations, including equations with constant coefficients, represent special cases of quasi-Bessel equations and, as just a corollary, we can derive the key points of the theory for fractional constant-coefficient equations.\\

The structure of this paper is as follows. \\
In Section \ref{SectBeyondFracBesEq} we present the fractional series solution for equation \eqref{QuasiBessel} and derive the steps to obtain it numerically.\\

In Section \ref{sectNecessityP_1} we cite the existence theorem from \cite{2020DuSl-3}, and prove our main existence result (Theorem \ref{MainThrm}), including the necessary condition claiming that the power of $x$ must match the highest fractional derivative. \\

Motivated by \cite{Rodrigues}, Section \ref{sectUniq} presents the uniqueness theorem for the initial value problem of equation \eqref{QuasiBessel}.\\

Section \ref{SectEqnConstCoeffs} covers the use of the presented methodology to identify a series solution for homogeneous fractional equations with constant coefficients and its further extension into equations with power function factors.  
As a particular case of the constructed theory we have another derivation of the existence results for fractional differential equations with constant coefficients with multiple derivatives of any order.\\

Section \ref{sectExamples} provides a few numerical examples which support the constructed theory. 
Numerical method which utilizes substitution \cite{DS2020DuSl-2} is used to cross check the identified series solution in Example 1 and extend it further. Examples 2-4 confirm that our methodology works on previously analytically solved problems with one fractional derivative presented in monograph by Kilbas, Srivastava, and Trujillo \cite{Kilbas}.\\

 Throughout this paper notation $D^\alpha u(x)$ is used for both, Riemann-Liouville and Caputo, types of fractional derivatives, whereas $D_C^\alpha u(x)$ refers to the Caputo derivative, and $D_R^\alpha u(x)$ is the Riemann-Liouville derivative.     
%

\section{Construction of the fractional series solution}\label{SectBeyondFracBesEq}
We go beyond equation \eqref{GenBesselEqn} and consider equation \eqref{QuasiBessel} in the following form:
\begin{equation}\label{GenBesselEqnBeyond-2}
	d_1 x^{\alpha_1}D^{\alpha_1}u(x)+\sum_{i=2}^{m}d_i x^{\alpha_i+p_i}D^{\alpha_i} u(x) + (x^\beta - \nu^2)u(x)=0. 
\end{equation}
 It is worth mentioning its particular cases:
\begin{eqnarray}
	\hspace*{-8mm} &a)& d_1 x^{\alpha_1}D^{\alpha_1}u(x)+\sum_{i=2}^{m}d_i x^{\alpha_i+p_i}D^{\alpha_i} u(x) + x^\beta u(x)=0, \ \nu =0, \nonumber \\
	\hspace*{-8mm} &b)& d_1 x^{\alpha_1}D^{\alpha_1}u(x)+\sum_{i=2}^{m}d_i x^{\alpha_i+p_i}D^{\alpha_i} u(x) +  u(x)=0, \
	\beta=\nu=0. \nonumber 
\end{eqnarray}
\begin{df}
Differential equation (\ref{GenBesselEqnBeyond-2}) is called a quasi-Bessel equation
provided that $\alpha_1 = \max\{\alpha_i\}$ and for all $i,\, 1\leq i\leq m$,
$d_i\in \mathbb{R}$, $\alpha_i \in \mathbb{R}^+=[0,\infty)$, $\beta\in \mathbb{Q^+}=\mathbb{R}^+\cap \mathbb{Q}$.
Shifting indices $p_i$, $2\leq i\leq m$, 
should be non-negative and almost-rational, i.e., for a fixed number $r\in \mathbb{R^+}$ 
the deviations $p_i$ must belong to the set 
\[
p_i\in \mathbb{Q}^+_r=r\cdot\mathbb{Q^+}=\{x:\ x=rq,\ q\in \mathbb{Q^+} \}.
\] \qed
\end{df}
This section is applicable to equations with both Caputo and Riemann-Liouville derivatives.  The only difference is the condition on acceptable $\gamma$ in \eqref{CondForC0} to generate a true solution.  For Riemann-Liouville it is $\gamma > -1$ and for Caputo -- $\gamma > \lceil \alpha_1 \rceil -1$ to assure existence of derivatives.\\

{\scshape Notation.} 
Let $n_{\max}$ be the integer ceiling for the highest non-integer derivative, 
i.e., $n_{\max}=\max{\{n_i\}}$ and $n_{\min}$ be the integer ceiling for the lowest non-integer derivative, 
i.e., $n_{\min}=\min{\{n_i\}}$ where \\ $n_i-1 <\alpha_i < n_i, n_i \in \mathbb{N}, i=1,...,m$. 

\begin{lm}\label{lemma2} (see \cite{2020DuSl-3})
	{\it For the existence of the series solution for equation \eqref{GenBesselEqnBeyond-2}  with Caputo derivatives, 
	it is necessary that solutions of the characteristic equation \eqref{CondForC0},  
	$\gamma > n_{\max}-1 \ge 0$. } \qed
\end{lm}
In the multi-term fractional Bessel equation \eqref{GenBesselEqn} it was possible to use $s=\beta$ as a step (increase) in the series solution because it allowed recursively express all $c_n$ \eqref{Eqnforc_n}.  In the more general case of equation \eqref{GenBesselEqnBeyond-2}, it is not possible since some coefficients $c_n$ do not get 'balanced out' and need to be forced to be zero, which leads to a trivial solution.  Therefore, we need to find the largest possible step $s \le \beta$ which allows us to balance out (express recursively) all coefficients. 
Thus, we proceed as follows.\\

Let us assume that the solution for equation \eqref{GenBesselEqnBeyond-2} is 
\begin{equation}\label{SolnFnB}
u(x)=\sum_{n=0}^{\infty}c_n x^{\gamma+sn},
\end{equation}
where $\gamma$ satisfies the equation
\begin{equation}\label{CondForC0B}
G(\gamma)=\sum_{i=1}^{m_1}\frac{d_i\Gamma(1+\gamma)}{\Gamma(1+\gamma-\alpha_i)}-\nu^2=\sum_{i=1}^{m_1}d_i Q(0,\alpha_i)-\nu^2=0,
\end{equation}
where $p_i=0, i=2,...,m_1$ and
\begin{equation}\label{DefQ}
Q(r,p)=\frac{\Gamma(1+\gamma+r)}{\Gamma(1+\gamma+r-p)}.
\end{equation}
It is clear that all roots of equation \eqref{CondForC0B} belong to a bounded interval. \\
In addition to $p_1=0$, several other terms could also have $p_i=0$. Let us call them the pure Bessel terms. For these terms the power of the factor $x^\alpha$ matches the order of the derivative $D^\alpha u(x)$.  Let $m_1$ be the number of pure Bessel terms in \eqref{GenBesselEqnBeyond-2}.  Then, all other terms have strictly positive shifted powers $p_i > 0$ for $i=m_1+1,...,m$.\\ 

By plugging expression \eqref{SolnFnB} into equation \eqref{GenBesselEqnBeyond-2}, we obtain 
\begin{eqnarray}
\sum_{n=0}^{\infty}\!\!\!\! &c_n&\!\!\!\! \left( \sum_{i=1}^{m_1}d_i x^{\gamma+s n}Q(ns,\alpha_i)+\sum_{i=m_1+1}^{m}d_i x^{\gamma+s n +p_i}Q(ns+p_i,\alpha_i) \right) \nonumber \\
&+&\!\!\! (x^\beta - \nu^2)\sum_{n=0}^{\infty}c_n x^{\gamma+s n}=0, \nonumber
\end{eqnarray}
or
\begin{eqnarray}\label{SeriesEq} 
\sum_{n=0}^{\infty} c_n x^{\gamma+s n}\!\! &\Bigg(\!\!&\ds\sum_{i=1}^{m_1} d_iQ(ns,\alpha_i)-\nu^2 \nonumber \\ 
&+&\!\!\ds\sum_{i=m_1+1}^{m} x^{p_i}d_iQ(ns+p_i,\alpha_i)+x^\beta \Bigg)=0. 
\end{eqnarray}
If the step $s$ is such that $\ds\frac{p_i}{s}=n_{p_i} \in \mathbb{N}$ and $\ds\frac{\beta}{s}=n_\beta \in \mathbb{N}$ then we obtain
\begin{eqnarray}\label{SeriesEqb}
\sum_{n=0}^{\infty} c_n x^{\gamma+s n} \Bigg( &\ds\sum_{i=1}^{m_1}&d_iQ(ns,\alpha_i)-\nu^2 \nonumber \\
+&\ds\sum_{i=m_1+1}^{m}&x^{s n_{p_i}}d_iQ(ns+p_i,\alpha_i)+x^{s n_\beta} \Bigg)=0. 
\end{eqnarray}\\

Step $s$ should be such that any powers of $x$ in quasi-Bessel equation \eqref{GenBesselEqnBeyond-2} are included in the set $\gamma+sn$.  We want to find the maximum possible step $s$ to minimize the number of steps and avoid having unnecessary zero coefficients $c_n$. Let us describe how step $s$ can be identified.  
\begin{itemize}
	\item From the definition of quasi-Bessel equation \eqref{GenBesselEqnBeyond-2} we have  $\beta \in \mathbb{Q}_r^+$, $p_i \in \mathbb{Q}_r^+$. Therefore we first convert them into rational numbers:\\
	$\ds \beta^0=\frac{\beta}{r}=\frac{a_1}{b_1}$; 
	$\ \ \ds p_i^0=\frac{0}{1}=\frac{a_i}{b_i},2 \le i \le m_1; \ \  p_i^0=\frac{p_i}{r}=\frac{a_i}{b_i}$, 
	$m_1 < i\le m$.  If $r\in \mathbb{Q^+}$, then $r=1$.
	\item Find the lowest common denominator $N_{LCD}=\text{LCD}\{b_i\}$,\\
	$m_1 < i\leq m$. 
	\item Calculate the acceptable step and corresponding shifts for $\beta$ and $p_i$:
	\begin{equation}\label{BetasNs0}
	s^0=\frac{1}{N_{LCD}}; \ n_\beta^0 = \frac{\beta^0}{s^0}\in\mathbb{N}; \  n_{p_i}^0=\frac{p_i^0}{s^0}\in\mathbb{N}, \ m_1<i\leq m.
	\end{equation}
	\item The identified parameters $\beta^0, p_i^0$, $m_1<i\leq m$ in \eqref{BetasNs0} can still have common factors.  To maximize step $s$ we need to identify their greatest common factor ($N_{gcf}$), adjust step $s$ and each parameter. Then, finally, we obtain: 
	\begin{equation}\label{BetasNs}
	s=s^0 \cdot N_{gcf};\ n_\beta = \frac{n_\beta^0}{N_{gcf}}; \ n_{p_i}=\frac{n_{p_i}^0}{N_{gcf}},\ m_1< i\leq m.
	\end{equation}
\end{itemize}
As an example, for equation
\[
2x^{2.4}D^{2.4} u(x) - 3x^{1.8} D^{1.5} u(x) +xD^{0.4} u(x) + (x^3-\nu^2) u(x)=0
\]
we have $d_1=2$, $d_2=-3$, $d_3=1$, $\alpha_1=2.4$, $\alpha_2=1.8$,
$\alpha_3=0.4$. Then $\beta=3$, $p_2=0.3$, $p_3=0.6$, and we obtain
$b_1=1$, $b_2=10$, $b_3=5$, their $N_{LCD}=10$. 
Thus,
\[
s^0=\frac1{N_{LCD}}=0.1,\ n_\beta^0 =\frac{\beta}{s^0}=30,\ n_{p_2}^0=\frac{p_2}{s^0}=3,\
n_{p_3}^0=\frac{p_3}{s^0}=6. 
\]
Since $N_{gsf}={\rm GCF}(30,3,6)=3$, then, finally, 
\[
s=s^0\cdot N_{gsf} =0.3,\ n_\beta=\frac{30}{3}=10,\ n_{p_2}=1,\ n_{p_3}= 2.
\]
Other examples are presented in Sections \ref{SectEqnConstCoeffs} and  \ref{sectExamples}.\\

Equation \eqref{SeriesEqb} can be re-written to identify more clearly the recursive relationship of the coefficients $c_n$:
\begin{eqnarray}\label{SeriesEqb2}
\hspace*{-5mm}&&\sum_{n=0}^{\infty} c_n x^{\gamma+s n} \left( \sum_{i=1}^{m_1}d_iQ(ns,\alpha_i)-\nu^2\right) + \nonumber \\
&&\!\!\!\!\sum_{i=m_1+1}^{m}\left(\sum_{n=n_{p_i}}^{\infty} c_n x^{\gamma+s n}d_iQ(ns+p_i,\alpha_i)\right)
+\sum_{n=n_\beta}^{\infty} c_n x^{\gamma+s n} =0. \: \: \:  \: \: \: 
\end{eqnarray}
As long as condition \eqref{CondForC0} for $x^\gamma$ at $n=0$ is satisfied, we can make $c_0$ equal to any number. For simplicity we assume $c_0=1$.  Then, to make equation \eqref{SeriesEqb2} valid, we need to offset $c_n$ by the sum of coefficients $c_n^{p_i}$ and $c_n^\beta$ 
\begin{equation*}
c_n=c_n^\beta+\sum_{i=m_1+1}^{m}c_n^{p_i}.
\end{equation*}
Coefficient $c_n$ should compensate previous like terms, the terms which are $n_\beta$ steps before $c_n$ together with the terms which are $n_{p_i}$ steps before $c_n$ from \eqref{SeriesEqb2}.\\  

These coefficients can be expressed as
 \begin{equation}\label{CoeffType1}
 c_n^{p_i}=-\frac{c_{n-n_{p_i}}\cdot d_{p_i} Q\left((n-n_{p_i})s,\alpha_{p_i}\right)} {\ds\sum_{i=1}^{m_1}d_iQ(ns,\alpha_i)-\nu^2},
 \end{equation}
 where $d_{p_i}$ and $\alpha_{p_i}$ correspond to $p_i, m_1 < i \le m$ from the original equation \eqref{GenBesselEqnBeyond-2} and
 \begin{equation}\label{CoeffType2}
 c_n^{\beta}=-\frac{c_{n-n_\beta}}{\ds\sum_{i=1}^{m_1}d_iQ(ns,\alpha_i)-\nu^2}.
 \end{equation}
These coefficients \eqref{CoeffType1} and \eqref{CoeffType2} are equal to zero up to the point $n_{p_i}$ for the terms $c_n^{p_i}$ and up to $n_\beta$ for $c_n^\beta$.  Hence, we combine (\ref{CoeffType1}), (\ref{CoeffType2}) and obtain the following formula for coefficients $c_n$: 
\begin{eqnarray}\label{FinalEqforCn}
c_n\!\!\!\! &=&\!\!U(n-n_\beta)c_n^{\beta}+\sum_{i=m_1+1}^{m}U(n-n_{p_i})c_n^{p_i} \nonumber \\
&=&\!\!-\frac{U(n-n_\beta)c_{n-n_\beta}+\ds\sum_{i=m_1+1}^{m}U(n-n_{p_i})c_{n-n_{p_i}} \cdot d_i Q((n-n_{p_i})s, \alpha_i) }{\ds\sum_{i=1}^{m_1}d_iQ(ns,\alpha_i)-\nu^2},\nonumber\\  
\end{eqnarray}
where $
U(n-k)=
\begin{cases}
1 \text { if } n\ge k \\
0 \text{ if } n < k
\end{cases}
$
and $Q(r,p), s, n_{p_i}, n_\beta$ are defined in \eqref{DefQ}, \eqref{BetasNs}.
This algorithm provides the fractional series solution for equation \eqref{GenBesselEqnBeyond-2}.

\section{Existence theorem}\label{sectNecessityP_1}
In this section we prove that the constructed fractional series solution 
and matching 
power $x^\alpha$ and the highest derivative (i.e., $p_1=0$) is the necessary and sufficient condition  to obtain the series solution \eqref{SolnFnB} where coefficients are defined in  \eqref{FinalEqforCn}.
Based on the formulas in the previous section it looks like the described logic provides a solution to equation \eqref{GenBesselEqnBeyond-2} if at least one $p_i=0, i=1,...,m$ (not necessarily $p_1=0$). However, this is not the case as shown in Theorem \ref{MainThrm} below. \\
	
Let $n_{m_0}=\ds\max_{1\le i \le m_0}{\{n_i\}}$ where $n_i-1 <\alpha_i < n_i, \ n_i \in \mathbb{N},  i=1,...,m_0$. 
Here $m_0 \le m_1$ is the number of fractional derivatives and $m_1-m_0$ is the number of integer derivatives with $p_i=0$ in equation \eqref{GenBesselEqnBeyond-2} or the number of integer derivatives in multi-term fractional Bessel equation \eqref{GenBesselEqn}.

\begin{thm}\label{Thrm1} \cite{2020DuSl-3}. 
	{\it There exists a series solution \eqref{SolnFn},\eqref{Eqnforc_n} for fractional equation \eqref{GenBesselEqn} with Caputo derivatives with all $d_i>0$ in any domain 
	$x \in [0,b], b \in \mathbb{R_+}$ if $\nu$ satisfies the inequality
	\begin{equation}\label{TrueNuCondinTheorem}
		\nu^2=\nu^2_{\min} \ge \Gamma(n_{m_0})\sum_{i=1}^{m_1}\frac{d_i}{\Gamma(n_{\max}-\alpha_i)}.
	\end{equation}
	%
	If $\nu = 0$ and $\beta \ge n_{m_0}-1$ then at least one solution in the form of series can always be found.\\
	If $\ds\alpha_{\max}^*=\max_{1\leq i\leq m_1}\{\alpha_i\}$ (the highest value of $\alpha_i$ with $p_i=0$) is fractional and $\nu_{\min} > 0$, then the series solution is unique.\\
	If $\alpha_{\max}^*$ is integer and $\nu_{\min} > 0$, then the series solution is unique provided that
	$n_{m_0} \ge \alpha_{\max}^*-1$. \\
	If $n_{m_0} < \alpha_{\max}^*-1$, then there may be several series solutions.} \qed
\end{thm}	

	\begin{rem} \label{RemarkRootNotSoln}
			If there exists $n$ such that $\gamma+sn$ is another root of \eqref{CondForC0} with step $s$ defined in \eqref{BetasNs}, then $\gamma$ does not generate solution \eqref{SolnFnB} for equation \eqref{GenBesselEqnBeyond-2} because in this case the series is divergent.\\
			
	It usually happens when $\nu=0$ and only $p_1=0$. Then the difference between $\gamma$ roots is exactly one.  If step $s$ is a fraction of one, smaller $\gamma$ root at some step falls onto a bigger $\gamma$ root and the series blows up.\\
	 
	For example, equation
		\begin{equation*}
			x^{1.5}D^{1.5}u(x)+x^{0.7}D^{0.5}u(x)+(x^{1.2}-0)u(x)=0
		\end{equation*}
		has characteristic equation
		\begin{equation*}
			\frac{\Gamma(1+\gamma)}{\Gamma(1+\gamma-1.5)}=0.
		\end{equation*}
		Hence, $\ds\beta=\frac{6}{5}, p_2=\frac{1}{5}, r=1$.  $N_{LCD}=5, s=\ds\frac{1}{5}, n_\beta=6, n_{p_2}=1$.  Roots are $\gamma_1=-0.5, \gamma_2=0.5$.  Therefore $\gamma_2=\gamma_1+5s$, which means that for $n=5$
		\begin{equation*}
			Q(5s,\alpha_1)=\frac{\Gamma(1+\gamma_1+5s)}{\Gamma(1+\gamma_1+5s-\alpha_1)}=\frac{\Gamma(1-0.5+1)}{\Gamma(1-0.5+1-1.5)}=\frac{\Gamma(-0.5)}{\Gamma(0)}=0,
		\end{equation*}
		which is in the denominator of $c_5$ in \eqref{FinalEqforCn} and makes $c_5=\infty$, therefore
		 $\gamma_1$ does not generate solution in the form of proposed series.\\
		 			
	Equation  \eqref{NoSolConst} in  Section \ref{SectEqnConstCoeffs} is also an example of this situation. \qed		
		\end{rem} 
Thus, if the characteristic equation \eqref{CondForC0} has several roots, then all roots except for the largest root need to be checked for validity.  
%
\begin{lm}\label{Lemma1} (Watson's lemma as a corollary of Sterling approximation for $\Gamma$-function).\\ 
	
	{\it For any $\alpha \in \mathbb{C}$,
	\begin{equation}\label{EulerGammasLim}
		\lim_{n\to\infty}\frac{\Gamma(n)n^\alpha}{\Gamma(n+\alpha)}=1.
	\end{equation}} \qed
	%
\end{lm}
%
\begin{thm}\label{MainThrm} 
	{\it Series solution \eqref{SolnFnB} with coefficients \eqref{FinalEqforCn} of fractional quasi-Bessel equation \eqref{GenBesselEqnBeyond-2} with $p_1=0$, $d_i>0$, $1\leq i\leq m_1$, converges and represents the solution of equation \eqref{GenBesselEqnBeyond-2} provided that the threshold  condition \eqref{TrueNuCondinTheorem} for $\nu$ is satisfied in the equations with Caputo derivatives.  No such threshold condition is required for the equations with Riemann-Liouville derivatives.  \\
	 
	 If $p_1 > 0$ but for some $i>1$ there exists at least one $p_i=0$, then the series diverges and a series solution in form \eqref{SolnFnB} does not exist for equations with both Caputo and Riemann-Liouville derivatives. }
\end{thm} 
\begin{proof}
 We evaluate this series starting from the element \\
   $c_n>\ds\max_{i}\{n_{p_i},n_\beta\}$, 
 we get:
	\begin{eqnarray*}
	|c_n| &\le&\left |\frac{c_{n-n_\beta}+\ds\sum_{i=m_1+1}^{m}c_{n-n_{p_i}} \cdot d_i Q((n-n_{p_i})s,\alpha_i) }{\ds\sum_{j=1}^{m_1}d_jQ(ns,\alpha_j)-\nu^2}
	   \right |  
	\end{eqnarray*}
   	\begin{eqnarray}
	   &\le& \left | \frac{c_{n-n_\beta}}{\ds\sum_{j=1}^{m_1}d_jQ(ns,\alpha_j)-\nu^2} \right| + \sum_{i=m_1+1}^{m}|d_i|\left |\frac{c_{n-n_{p_i}}\cdot Q((n-n_{p_i})s,\alpha_i) } {\ds\sum_{j=1}^{m_1}d_jQ(ns,\alpha_j)-\nu^2}\right | \nonumber \\
	   &\le& \left | \ds\prod_{k=0}^{n-n_\beta} \frac{1} {\ds\sum_{j=1}^{m_1}d_iQ((k+n_\beta)s,\alpha_j)-\nu^2} \right |  \label{TrmFstTerm} \\
	   &+&\sum_{i=m_1+1}^{m}|d_i|^{n-n_{p_i}}\left |
 	\ds\prod_{k=0}^{n-n_{p_i}}\frac{Q(ks, \alpha_i)} {\ds\sum_{j=1}^{m_1}d_jQ((k+n_{p_j})s,\alpha_j)-\nu^2} \right |. \label{TrmSndTerm}
	\end{eqnarray}
	 Let us evaluate the first term \eqref{TrmFstTerm}. Keeping in mind Remark \ref{RemarkRootNotSoln}, the denominator is not equal to zero except for the roots of the characteristic equation \eqref{CondForC0} and therefore since $d_j > 0$ for $1 \le j \le m_1$ and infinite increase of $Q(k,\alpha)$, we can ignore $\nu^2=\const$  because it does not effect convergence.  Based on expression \eqref{DefQ} for $Q(r,p)$ and Euler's formula \eqref{EulerGammasLim} we can conclude 
	\begin{eqnarray*}
		&&\left|\frac{1}{\ds\sum_{j=1}^{m_1}d_jQ((k+n_\beta)s,\alpha_j)}\right|\le\left|\frac{1}{d_1Q((k+n_\beta)s,\alpha_1)}\right| \\
		&&\ \ =\left|\frac{\Gamma(1+\gamma+(k+n_\beta)s-\alpha_1)}{d_1\Gamma(1+\gamma+(k+n_\beta)s)}\right|=O\left(\left|\frac{1}{d_1(1+\gamma+(k+n_\beta)s)^{\alpha_1}}\right|\right) \\
		&& \ \ \le O\left(\left|\frac{1}{d_1((k+n_\beta)s)^{\alpha_1}}\right|\right), \ k \to \infty \ \text{ (since } 1+\gamma = \const > 0).  
	\end{eqnarray*}
	Hence we derive  $(d_j > 0, \ 1\le j\le m_1; \ \ \alpha_1 > \alpha_i, \ 2\le i\le m)$:
 	\begin{eqnarray}\label{Trm2FstTerm}
		&&\left|\ds\prod_{k=0}^{n-n_\beta}\frac{1}{\ds\sum_{j=1}^{m_1}(d_iQ((k+n_\beta)s,\alpha_j)-\nu^2)}\right| \nonumber \\
		&&\ \ \ \le O\left(\left[\frac{1}{d_1^{n-n_\beta+1}\big(n_\beta s\big)\big((1+n_\beta) s\big)\big((2+n_\beta)s\big)\cdot ...\cdot\big((n-n_\beta+n_\beta)s\big)}\right]^{\alpha_1}\right) \nonumber \\
		&&\ \ \ =O\left(\left[\frac{(n_\beta-1)!}{(d_1 s)^{n-n_\beta+1}n!}\right]^{\alpha_1}\right)
	\end{eqnarray}
	Evaluation of the second term \eqref{TrmSndTerm} is similar:
	\begin{eqnarray*}
		&&|Q(ks,\alpha_i)|=\left|\frac{\Gamma(1+\gamma+ks)}{\Gamma(1+\gamma+ks-\alpha_i)}\right|=O((1+\gamma+ks)^{\alpha_i}).\\
	\end{eqnarray*}
	Since the largest root $\gamma$ is less than some constant, then we can say that \\ $0 < 1+\gamma \le c_\gamma s =  \const, \ c_\gamma\in \mathbb{N}$.	Therefore,
	\begin{eqnarray}\label{Thm2NumSecTerm}
	    &&\left|\ds\prod_{k=0}^{n-n_{p_i}}Q(ks, \alpha_i)\right| =O\big(\big[(1+\gamma)(1+\gamma+s)(1+\gamma+2s)\cdot ... \nonumber \\
	    &&\ \ \ \ \ \ \ \ \ \ \ \ \ \ \ \ \ \ \ \ \ \ \ \ \ \ \ \ \cdot (1+\gamma+(n-n_{p_i})s)\big]^{\alpha_i}\big) \nonumber \\
	    &&\le O\big(\big[c_\gamma s(c_\gamma s+s)(c_\gamma s+2s)\cdot ... \cdot (c_\gamma s+(n-n_{p_i})s)\big]^{\alpha_i}\big) \nonumber  \\
	    &&=O\left(\left[s^{n-n_{p_i}+1}(c_\gamma+1)(c_\gamma+2)\cdot ... \cdot ((c_\gamma+n-n_{p_i}))\right]^{\alpha_i}\right) \nonumber \\
		&&=O\left(\left[s^{n-n_{p_i}+1}\frac{(c_\gamma+n-n_{p_i})!}{c_\gamma!}\right]^{\alpha_i}\right). 
	\end{eqnarray}	
	Similarly to \eqref{Trm2FstTerm}, we obtain
	\begin{eqnarray*}
		\left|\ds\prod_{k=0}^{n-n_{p_j}}\frac{1}{\ds\sum_{j=1}^{m_1}(d_jQ((k+n_{p_j})s,\alpha_j)-\nu^2)}\right|
		=O\left(\frac{(n_{p_j}-1)!}{(d_1 s)^{n-n_{p_j}+1}n!}\right)^{\alpha_1}.
	\end{eqnarray*}
	Given that $n_\beta, n_{p_i}, d_i^{n_\beta}, s^{n_\beta}, d_i^{n_{p_i}}, s^{n_{p_i}}, m_1 < i \le m,$ and $c_\gamma$ are constants, we can ignore them in the evaluation of convergence.\\
	
	Expression \eqref{Thm2NumSecTerm} becomes
	\begin{eqnarray}
		&&O\left(\left[s^n (c_\gamma+n-n_{p_i})!\right]^{\alpha_i}\right)=O\left(\left[s^n (n-n_{p_i})!P_{c_\gamma}\right]^{\alpha_i}\right), \ \text{ where } \nonumber \\
		&& P_{c_\gamma}=(n-n_{p_i}+1)(n-n_{p_i}+2)...(n-n_{p_i}+c_\gamma) \text{ is a polynomial of degree $c_\gamma$}.    \nonumber \\
		&&\text{Hence, } \nonumber \\
		&& O\left(\left[s^n (c_\gamma+n-n_{p_i})!\right]^{\alpha_i}\right)=O\left(\left[s^n (n-n_{p_i})!n^{c_\gamma}\right]^{\alpha_i}\right).
	\end{eqnarray}
	Consequently, we obtain  	
  	\begin{eqnarray}
  	|c_n| &=& O\left(\frac{1}{((d_1 s)^{n}n!)^{\alpha_1}}\right) +O\left(\sum_{i=m_1+1}^{m}|d_i|^n
  	\frac{[s^{n}(n-n_{p_i})!n^{c_\gamma}]^{\alpha_i}}{(d_1^n s^n n!)^{\alpha_1}}\right) \nonumber \\
  	&\le& O\left(\frac{1}{(d_1^n s^n n!)^{\alpha_1}}\right) +O\left(\sum_{i=m_1+1}^{m}
  	\frac{|d_i|^n n^{c_\gamma \alpha_i}}{d_1^{n\alpha_1}(s^n n!)^{\alpha_1-\alpha_i}}\right).
   \label{LimProdK} 
	\end{eqnarray}
	Since $\alpha_i < \alpha_1, \ m_1 < i \le m$, then we can conclude that 
	\begin{equation}
	|c_n|x^{\gamma+ s n} \le O\left(\frac{x^{sn}}{(d_1^n s^n n!)^{\alpha_1}}\right)+
	 O\left(\sum_{i=m_1+1}^{m}\frac{(x^s |d_i|)^n n^{c_\gamma \alpha_i}}{d_1^{n\alpha_1}(s^n n!)^{\alpha_1-\alpha_i}}\right), 
	\end{equation}
	Finally, series solution \eqref{SolnFnB} for equation \eqref{GenBesselEqnBeyond-2} uniformly converges on each bounded interval $x\in[\epsilon, b], \epsilon > 0$, and the solution is differentiable up to any order for $x > 0$.\\

	Conversely, if $p_1>0$ but $p_i=0, i=2,...,m_1$ then in this case  
	\begin{eqnarray}
		c_n&=&U(n-n_\beta)c_n^{\beta}+\sum_{i=m_1+1}^{m}U(n-n_{p_i})c_n^{p_i}+U(n-n_{p_1})c_n^{p_1} \nonumber \\
		&=&-\frac{U(n-n_\beta)c_{n-n_\beta}+\ds\sum_{i=m_1+1}^{m}U(n-n_{p_i})c_{n-n_{p_i}} \cdot d_i Q((n-n_{p_i})s, \alpha_i)}{\ds\sum_{i=2}^{m_1}d_iQ(ns,\alpha_i)-\nu^2} \nonumber \\
			&&-\frac{U(n-n_{p_1})c_{n-n_{p_1}} \cdot Q((n-n_{p_1})s, \alpha_1)} {\ds\sum_{i=2}^{m_1}d_iQ(ns,\alpha_i)-\nu^2}. \nonumber 
	\end{eqnarray}
	Per expression \eqref{DefQ}  and Euler's formula \eqref{EulerGammasLim} 
	and assuming that\\ $\alpha_m=\ds\max_{2 \le i \le m_1}\{\alpha_i\}$ and $d_m$ corresponds to $\alpha_m$, we get as in the first part of the theorem 
 	\begin{eqnarray}\label{Thm2DivergPart}
  	|c_n| &=& O\left(\frac{1}{(d_m^ns^{n}n!)^{\alpha_m}}\right) +O\left(\sum_{i=m_1+1}^{m}
	\frac{d_i^{n}(s^n(n-n_{p_i})!)^{\alpha_i}n^{ c_\gamma \alpha_i}}{d_m^{n\alpha_m} (s^n n!)^{\alpha_m}}\right) \nonumber  \\
 	&+& O\left(\frac{d_1^{n}n^{ c_\gamma \alpha_1}}{d_m^{n\alpha_m}(s^n n!)^{\alpha_m-\alpha_1}}\right).
	\end{eqnarray}
	Since $\alpha_1 > \alpha_m$ then the last addend dominates two other terms in \eqref{Thm2DivergPart} and is divergent.
	Hence, $|c_n| \to \infty$ as $n\to\infty$ and the series solution is divergent for any value of $x > 0$ in \eqref{SolnFnB}. \\

	Thus, condition $p_1=0$ is necessary for the existence of a series solution. As an additional remark, we should point out that the positivity of coefficients $d_i$ is required only for the terms where the order of the derivative is equal to the power of multiplier $x$. In this case characteristic equation (\ref{CondForC0}) is guaranteed to have real roots.
\end{proof}

\begin{rem} 
	For equations with Caputo derivative based on conditions in Theorem \ref{Thrm1}:
	\begin{itemize}
		\item If $p_1=0, \nu > 0$ and $\alpha_1$ is fractional then the found series solution is unique up to a constant,
		\item It $p_1=0$ and $\alpha_1$ is integer equation \eqref{GenBesselEqnBeyond-2} may have multiple solutions.
	\end{itemize}\qed
\end{rem}

\begin{rem} 
	Root $\gamma$ in the solution, which is calculated in \eqref{CondForC0B},  depends solely on the terms in equation \eqref{GenBesselEqnBeyond-2} with $p_i = 0$. \qed
\end{rem}
\begin{rem} 
If in addition to $p_1=0$ any other $p_i=0$, it only improves convergence since the denominator in the calculations of the coefficients \eqref{FinalEqforCn} has more terms and $c_n$ faster converges to zero. \qed 
\end{rem}

\section{Uniqueness for the initial value problem}\label{sectUniq} 
We can expand Theorem 2 in \cite{2020DuSl-3} for multi-term fractional Bessel equation  to  the fractional quasi-Bessel equation \eqref{GenBesselEqnBeyond-2} with positive coefficients $d_i$.

\begin{thm}{}\label{Thrm2} 
	{\it The initial value problem for fractional equation \eqref{GenBesselEqnBeyond-2} with Caputo derivatives with the domain $x \in [0,b]$  and initial conditions $u^{(j)}(0)=u^{(j)}_0, j=0,1,...,l-1,$ where $\ds l=\lceil \alpha_1 \rceil $, has a unique solution for every $\nu \in \mathbb{R}$ such that
	\begin{equation}\label{NuCondinTheorem}
	\nu^2 > b_1^{\beta}+\sum_{i=1}^{m}q_i |d_i| b_1^{n_i+p_i}
	\end{equation}
	where
	\begin{eqnarray}
	b_1=\max\{1,b\} \text { and } 
	q_i=
	\begin{cases} 
	\ds \frac{1}{\Gamma(n_i-\alpha_i)(n_i-\alpha_i+1)} & ,\alpha_i < n_i \\
	1 & ,\alpha_i= n_i  
	\end{cases}.
	\end{eqnarray}}
\end{thm}
\begin{proof}
	We denote by $C^l$ the space of functions $u$ which are $l$ times continuously differentiable on 
	$[0,b]$ with the norm 
	\begin{equation*}
	||u||_{C^l} = \sum_{k=0}^{l}||u^{(k)}||_C = \sum_{k=0}^{l}\max_{x\in [0,b]} |u^{(k)}(x)|.  
	\end{equation*}  
	Operator $T:C^l \mapsto C$ is defined as follows
	\begin{equation*}
	(Tu)(x)=\frac{1}{\nu^2}\left[x^\beta u(x) + \sum_{i=1}^{m}d_i x^{\alpha_i+p_i}D_C^{\alpha_i} u(x)\right].
	\end{equation*}
	Then, we can write general equation in the form $u(x)=(Tu)(x)$.  Taking into account formulas 2.4.24-26  in Corollary 2.3 for Caputo derivative in Chapter 2 from \cite{Kilbas} for the interval $[0,b]$, we have: 
	\begin{equation*}
	||D_C^{\alpha_i} u||_C \le k_{\alpha_i} ||u||_{C^l}, \text{ where } k_{\alpha_i} = \frac{b ^{n_i-\alpha_i}}{\Gamma(n_i-\alpha_i)(n_i-\alpha_i+1)}, \text { when } \alpha_i < n_i
	\end{equation*}
	and
	\begin{equation*}
	||D_C^{\alpha_i}u||_C=||u||_{C^l}, \text { when } \alpha_i=n_i.
	\end{equation*}
	Therefore,
	\begin{equation*}
	||D_C^{\alpha_i}u||_C \le q_i b^{n_i-\alpha_i}||u||_{C^l},
	\end{equation*}
	and we obtain
	\begin{eqnarray}
	||Tu_1-Tu_2||_C &=& \frac{1}{\nu^2}||x^\beta (u_1(x)-u_2(x)) + \sum_{i=1}^{m}d_i x^{\alpha_i+p_i} D_C^{\alpha_i}(u_1-u_2)(x)||_C\nonumber \\
	&\le& \frac{1}{\nu^2}\left(b_1^\beta||u_1-u_2||_C +\sum_{i=1}^{m} |d_i| b_1^{\alpha_i+p_i} q_i b^{n_i-\alpha_i} ||u_1-u_2||_{C^l} \right) \nonumber \\
	&\le& \frac{1}{\nu^2}\left(b_1^\beta +\sum_{i=1}^{m} q_i |d_i| b_1^{n_i+p_i}\right)||u_1-u_2||_{C^l}.
	\end{eqnarray}
	Having condition \eqref{NuCondinTheorem} we conclude that $T$ is a contraction.  We apply the Banach fix point theorem to complete the proof.
\end{proof}

\section{Equations with constant coefficients}\label{SectEqnConstCoeffs}
Let us consider equation with constant coefficients at each derivative: 
\begin{equation}\label{EqnConstCoeff}
\sum_{i=1}^{m}d_i D^{\alpha_i}u(x) + u(x)=0, \alpha_1 > \alpha_i > 0, \ i=2,...,m
\end{equation}
and assume that there exists a real positive number $r$ such that for all $i$, $1\leq i\leq m$, 
$\alpha_i \in \mathbb{Q}_r^+ =r \mathbb{Q^+}$
We multiply each term by $x^{\alpha_1}$.  Then we obtain
\begin{equation}\label{EqnConstCoeff2}
d_1 x^{\alpha_1}D^{\alpha_1}u(x)+\sum_{i=2}^{m}d_i x^{\alpha_1}D^{\alpha_i}u(x) + (x^{\alpha_1}-0)u(x)=0, 
\end{equation}
which is the quasi-Bessel equation \eqref{GenBesselEqnBeyond-2} with $\nu=0, \beta=\alpha_1$. Let us apply the methodology from Section \ref{SectBeyondFracBesEq}. In this case characteristic equation \eqref{CondForC0} becomes
\begin{equation}\label{CondC0ConstCoeffs}
\frac{\Gamma(1+\gamma)}{\Gamma(1+\gamma-\alpha_1)}=0,
\end{equation}  
and we arrive at roots $\gamma = \alpha_1-k, k\geq 1$.  Since only one term has matching power of $x$, then in this case the roots do not depend on the coefficients $d_i$.\\

If $\alpha_1$ is integer then \eqref{CondC0ConstCoeffs} provides $\alpha_1$ roots. In the case when $\alpha_1$ is truly fractional, i.e. $\lceil \alpha_1 \rceil \ne \alpha_1$, we can see that all roots $\gamma < \lceil \alpha_1 \rceil$ and therefore characteristic equation does not produce any solutions for constant-coefficient equations with Caputo derivatives unless $\gamma$ is a non-negative integer.  	
In this case, the Caputo derivative $D_C^{\alpha_1} x^j$ is equal to zero and solutions can be found without the use of characteristic equation. 
For example, for the simple equation 
\begin{equation}\label{Example1KilbCap}
D^{\alpha}_C u(x) - \lambda u(x)=0, \ 0<\alpha<1, \lambda \in \mathbb{R}
\end{equation}
solution exists \cite{Kilbas} as a specific series -- the Mittag-Leffler function  
\begin{equation*}
u(x)=b E_\alpha [\lambda x^\alpha] = b \sum_{n=0}^{\infty}\frac{\lambda^n x^{n\alpha}}{\Gamma(1+\alpha n)}, b=u(0) \in \mathbb{R}.
\end{equation*}
Our formula for coefficients produces identical result with $\beta=\alpha$,\\ $d_1=-\frac{1}{\lambda}, c_0=b$,
\[
c_{n+1}
 =\frac{c_n\lambda\Gamma(1+\alpha n)}{\Gamma(1+\alpha (n+1)}=\frac{\lambda^{ n+1}}{\Gamma(1+\alpha(n+1))}.
\]\\

In the case of equation \eqref{EqnConstCoeff2} with Riemann-Liouville derivatives the restriction on $\gamma$ is much milder: $\gamma > -1$, and, therefore, each root\\ $\gamma = \alpha_1-k, k=1,...,n_{\max}$ generates a solution for equation \eqref{EqnConstCoeff2} as long as it does not represent a root that becomes another root at some step $n$ as described in Remark \ref{RemarkRootNotSoln}.  In this case only the solution that corresponds to the largest root $\gamma$ exists. \\

{\scshape Example}. Characteristic equation \eqref{CondC0ConstCoeffs} of the equation
\begin{eqnarray}
  &&d_1 D_R^{2.1}u(x) + d_2 D_R^{1.4}u(x) + d_3 D_R^{0.7}u(x) + u(x)=0, \text{ or } \nonumber \\
  &&d_1 x^{2.1}D_R^{2.1}u(x) + d_2 x^{1.4+0.7}D_R^{1.4}u(x) + d_3 x^{0.7+1.4}D_R^{0.7}u(x) \nonumber \\
  &&+(x^{2.1}-0)u(x)=0, \nonumber \\
  && \alpha_1=2.1, \ \ p_2=0.7, \ \ p_3=1.4, \ \ \beta = 2.1, \nonumber
\end{eqnarray} 
has three roots: $\gamma_1=-0.9, \gamma_2=0.1$ and $\gamma_3=1.1$. The orders of the derivatives yield step $s=0.7$. Therefore $\gamma_2 \ne \gamma_1+sn, \gamma_3 \ne \gamma_1+sn, \gamma_3 \ne \gamma_2+sn$ for any $n\in \mathbb{N}$. Consequently, we
can construct three different series solutions based on these roots
\begin{eqnarray}
	u_1(x)=\sum_{n=0}^{\infty}c_n x^{-0.9+sn} \label{ExSol1} \\
	u_2(x)=\sum_{n=0}^{\infty}c_n x^{0.1+sn} \label{ExSol2} \\
	u_3(x)=\sum_{n=0}^{\infty}c_n x^{1.1+sn}. \label{ExSol3} 
\end{eqnarray}\\

In the case of almost the same equation
\begin{equation*}
  d_1 D_R^{2.1}u(x) + d_2 D_R^{1.5}u(x) + d_3 D_R^{0.7}u(x) + u(x)=0, 
\end{equation*}
which turns into
\begin{equation}\label{NoSolConst}
  d_1 x^{2.1}D_R^{2.1}u(x) + d_2 x^{1.5+0.6}D_R^{1.5}u(x) + d_3 x^{0.7+1.4}D_R^{0.7}u(x)  
  + x^{2.1}u(x)=0,  
\end{equation}
we obtain the same three roots $\gamma_1=-0.9, \gamma_2=0.1$ and $\gamma_3=1.1$, but the step here becomes $s=0.1$. In this case, $\gamma_2=0.1=-0.9+10\cdot 0.1=\gamma_1+10s$, which makes $c_{10}=\infty$, which invalidates series solution \eqref{ExSol1};   $\gamma_3=1.1=0.1+10\cdot 0.1=\gamma_2+10s$, which also makes $c_{10}=\infty$, and therefore invalidates series solution \eqref{ExSol2}  (see detailed explanation in example in Remark \ref{RemarkRootNotSoln}).  Consequently, neither $\gamma_1$ nor $\gamma_2$ represent a root which can be used to generate a series solution in the proposed form. \\

Root  $\gamma_3=\alpha_1-1$ ($\gamma_3=1.1$ in our example) is the largest and, in view of Remark \ref{RemarkRootNotSoln}, cannot be eliminated and is always valid.  Therefore, in the case of equation \eqref{EqnConstCoeff} with Riemann-Liouville derivatives we can always generate at least one series solution 
\begin{equation}\label{SolnFnC}
u(x)=\sum_{n=0}^{\infty}c_n x^{\alpha_1-1+sn}
\end{equation}
since, unlike roots $\gamma_1, \gamma_2$ of the characteristic equation, root $\gamma_3$ does not generate blow-up of the series solution.\\

Step $s$ is described in Section \ref{SectBeyondFracBesEq} for our modified equation \eqref{EqnConstCoeff2}, in which \\ $p_i=\alpha_1-\alpha_i, \beta=\alpha_1, \nu =0$, and the coefficients are calculated based on equation \eqref{FinalEqforCn}. \qed \\

As we see, constant-coefficient equations with Rieman-Liouville derivatives can be transformed into  quasi-Bessel equations.  Examples 2 and 3 in Section \ref{sectExamples} provide additional demonstrations of this ansatz.\\

Besides constant coefficients, the methodology described above is also applicable to equations with some terms having factors as powers of $x$.  We assume that $\beta_1$, the power of $x$ at the highest derivative $\alpha_1$, satisfies  inequality $\beta_1 \le \alpha_1$.  We consider equation 
\begin{eqnarray}\label{EqnConstXpowers}
&&\sum_{i=1}^{m}d_i x^{\beta_i} D_R^{\alpha_i}u(x) + x^{\delta}u(x)=0, \\
&&\alpha_1 > \alpha_i > 0, i=2,...,m, \ \alpha_1 \ge \beta_1, \ \alpha_1-\beta_1\geq \alpha_i-\beta_i, \nonumber 
\end{eqnarray}
where $\alpha_i, \beta_i, \delta \in \mathbb{Q}_r^+ =r \mathbb{Q^+}$ for some $r\in\mathbb{R^+}$.
We multiply each term by $x^{\alpha_1-\beta_1}$ and obtain
\begin{equation}\label{EqnXpowers2}
d_1 x^{\alpha_1}D_R^{\alpha_1}u(x)+\sum_{i=2}^{m}d_i x^{\alpha_1-\beta_1+\beta_i}D_R^{\alpha_i}u(x) + (x^{\alpha_1-\beta_1+\delta}-0)u(x)=0, 
\end{equation}
which is a quasi-Bessel equation.
In this case  $\nu=0, \beta=\alpha_1-\beta_1+\delta$.  Then, we can apply the above described methodology for equations with constant coefficients to this equation. In this case characteristic equation \eqref{CondForC0} is exactly the same as for equation with constant coefficients \eqref{CondC0ConstCoeffs}, and the approach for solving equations with constant coefficients applies verbatim.\\

Examples of equations with constant coefficients and with terms having powers of $x$ as factors are in Section \ref{sectExamples}, Examples 2-4.

\section{Examples of some equations and their solutions.}\label{sectExamples}
{\scshape Example 1} (quasi-Bessel equation with Caputo derivatives).\\

Let us consider equation
\begin{equation}\label{Example1}
1.5x^{1.5}D_C^{1.5}u(x)-1.2x^{1.9}D_C^{1.1}u(x) + 3xD_C^{0.5}u(x)+(x^2-\nu^2)u(x)=0.
\end{equation}
Here $\beta=2, d_1=1.5, d_2=-1.2, d_3=3, \alpha_1=1.5, \alpha_2=1.1, \alpha_3=0.5, \\ p_2=0.8, p_3=0.5$.\\

\begin{figure}[H]
	\begin{center}
		\includegraphics[width=6cm]{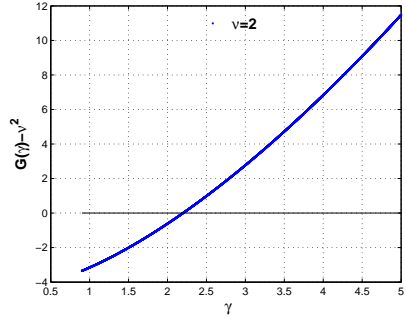}
		\caption{Function $G(\gamma)-\nu^2$ for equation \eqref{Example1} with $\nu=2$.}	\label{GyforEx1}
	\end{center}
\end{figure}

\begin{figure}[H]
	\begin{center}
		\includegraphics[width=6cm]{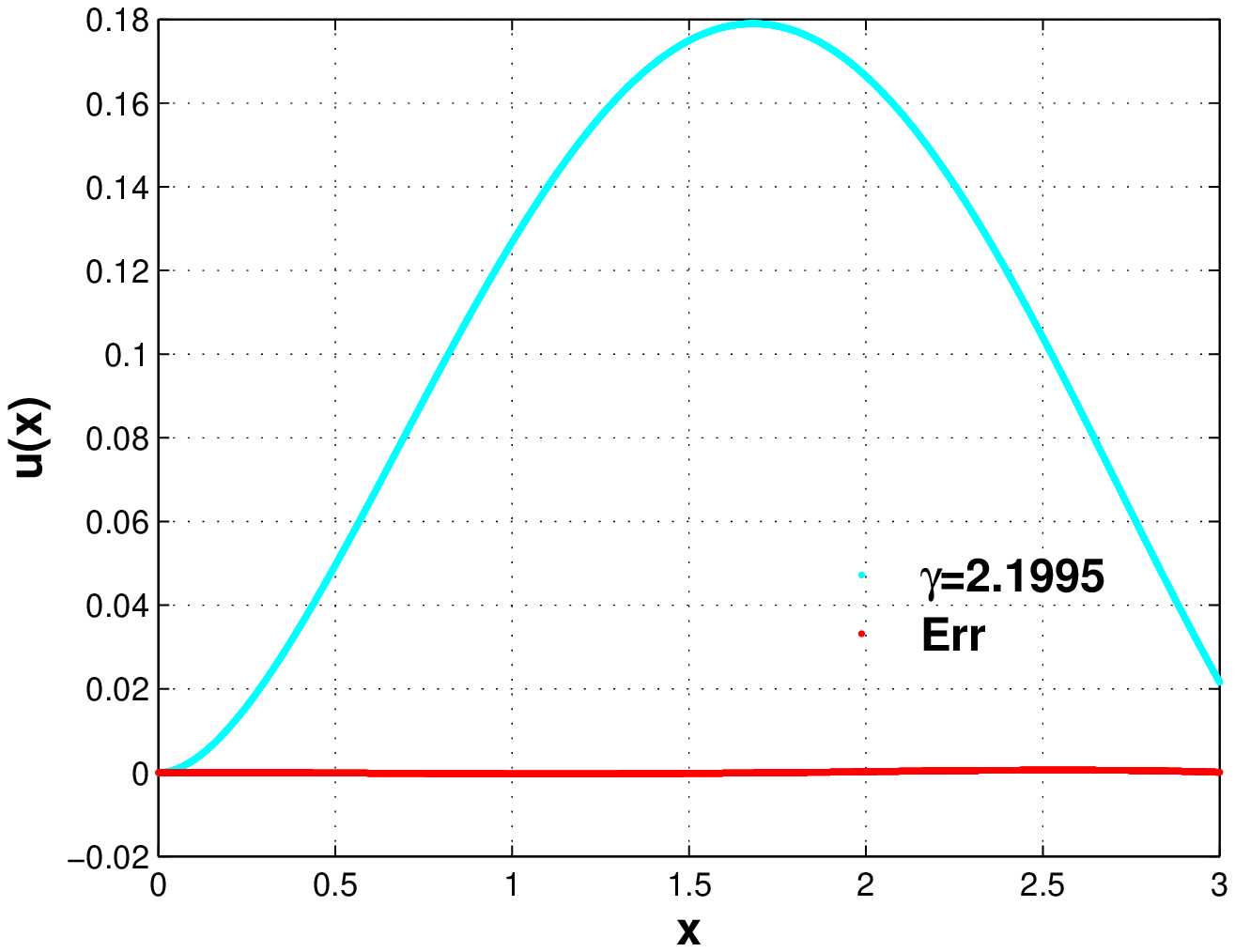}
		\includegraphics[width=6cm]{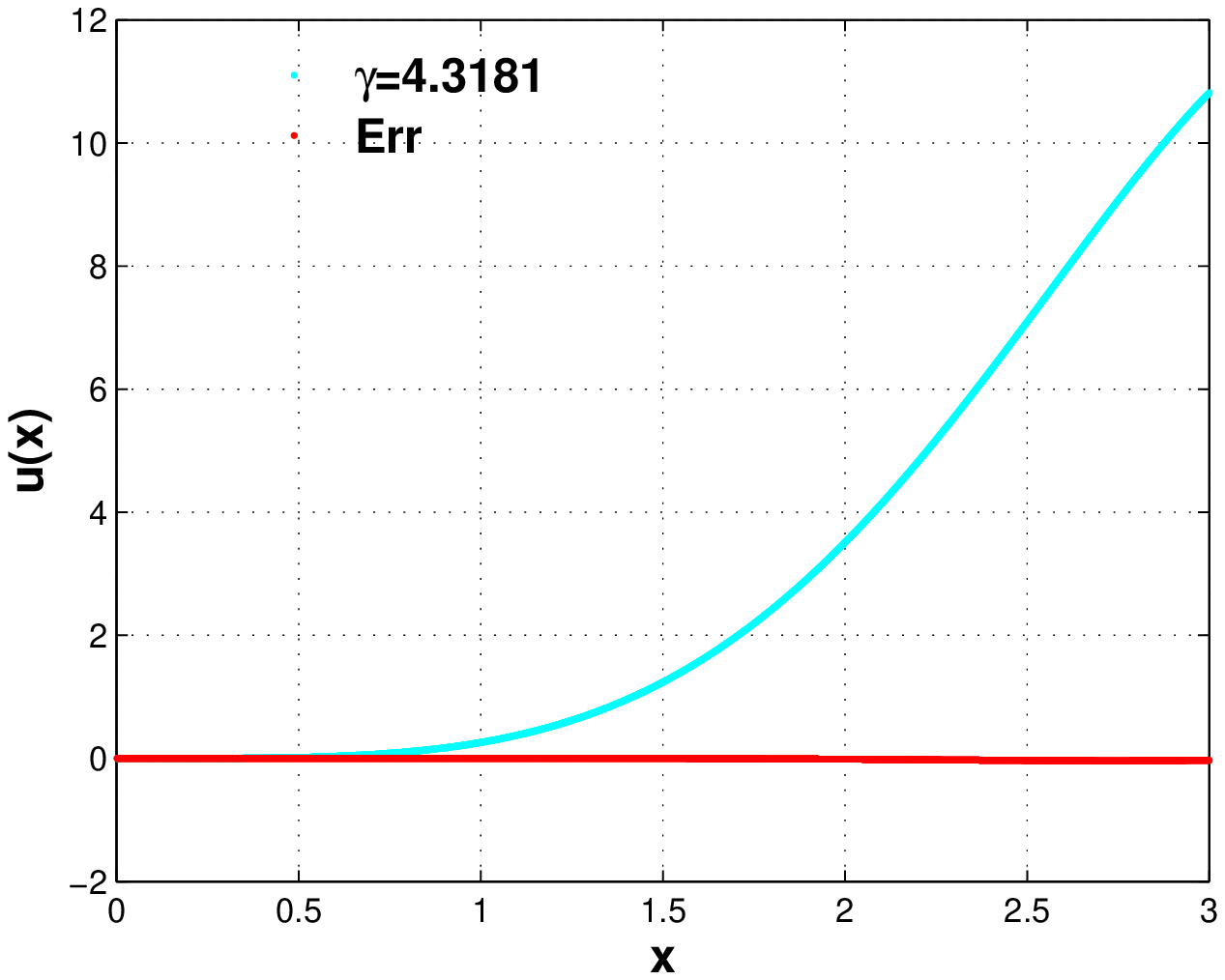}
\caption{Solution for equation in Example 1.  Red line close to zero is the check for the accuracy of the solution. Step $h=0.001$.}\label{SolEx1}	
	\end{center}
\end{figure}

Equation \eqref{CondForC0B} becomes: 
\begin{equation}\label{CondForC0Ex1}
G(\gamma)=\frac{1.5\Gamma(1+\gamma)}{\Gamma(1+\gamma-1.5)}-\nu^2=0.
\end{equation}
The graph of the expression on the left side of equation \eqref{CondForC0Ex1} is in Figure \ref{GyforEx1}. It is the same for any $\nu$ except for the $\nu^2$ shift down difference. To satisfy equation \eqref{CondForC0Ex1} for $\nu=2$ we found $\gamma=2.1995$, for $\nu=3.5$ -- $\gamma=4.3181$ (see \cite{2020DuSl-3} for the detail explanation on how to find all $\gamma$). \\
Other parameters involved in the process as described before are:
\begin{itemize}
	\item Since all $p_i, \beta \in \mathbb{Q^+}$ we get $\ds p_2^0=p_2=0.8=\frac{4}{5};$\\ 
	$\ds p_3^0=p_3=0.5=\frac{1}{2}; \ \ \ s^0 = \beta=2=\frac{2}{1}$. 
	\item The lowest common denominator $N_{LCD}=\text{LCM}\{5,2,1\}=10$.
	\item $\ds s=\frac{1}{N_{LCD}}=\frac{1}{10}=0.1, n_{p_2}=\frac{p_2}{s}=\frac{0.8}{0.1}=8, n_{p_3}=\frac{p_3}{s}=\frac{0.5}{0.1}=5$, \\$n_\beta = \ds\frac{\beta}{s}=\frac{2}{0.1}=20, (N_{gcf}=1)$.
\end{itemize} 
The solutions are represented in Figure \ref{SolEx1}.  Red line - recalculation of equation \eqref{Example1} by plugging in the calculated solution $u(x)$ into the equation.  Fractional derivatives are calculated using the substitution method for calculation of Caputo fractional derivatives as in \cite{DS2020DuSl-2}.
It is important to understand that the closer $\nu$ is to the minimum acceptable $\nu$ as explained in \cite{2020DuSl-3}, the less accurate the result is due to loss of accuracy in the calculation of fractional derivative.  \\

\begin{figure}[H]
	\begin{center}
		\includegraphics[width=6cm]{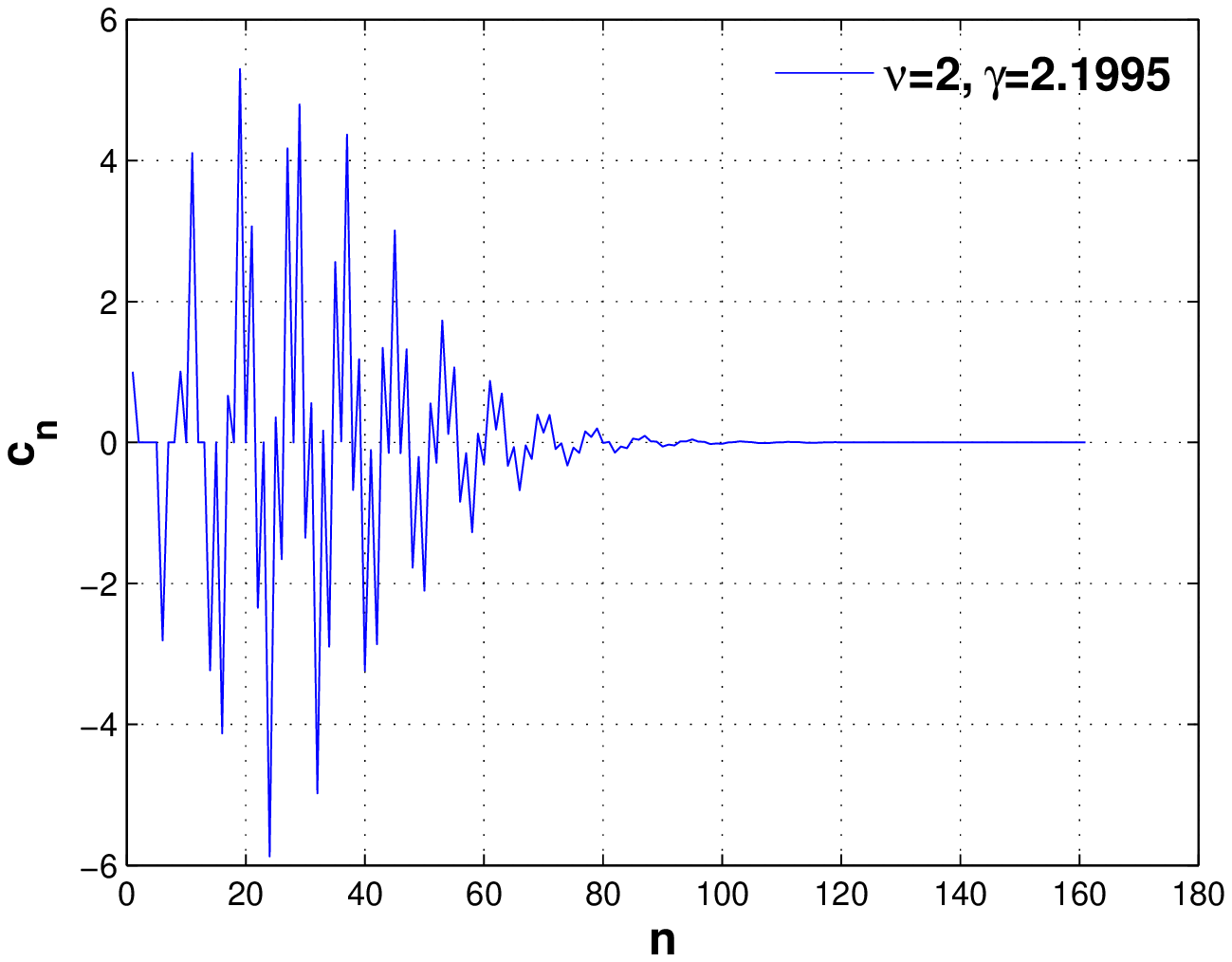}
		\includegraphics[width=6cm]{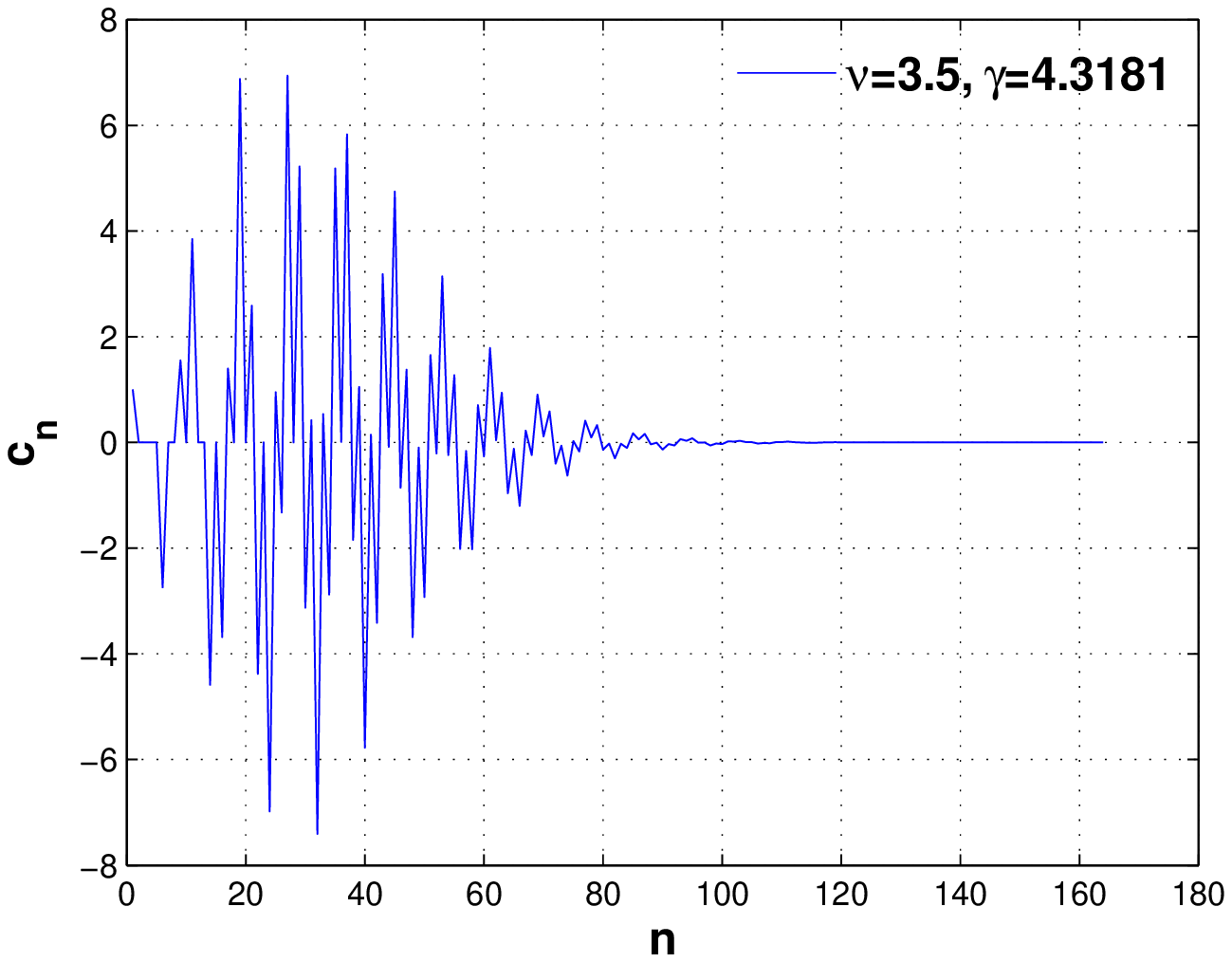}
		\caption{Example 1. Distribution of coefficients ($c_n > 10^{-5}$).}	\label{GraphC0Ex1}
	\end{center}
\end{figure}

Figure \ref{GraphC0Ex1} reflects distributions of $c_n$ for two cases of $\nu$.  It's easy to see that coefficients converge to zero very fast.

\begin{figure}[H]
	\begin{center}
		\includegraphics[width=6cm]{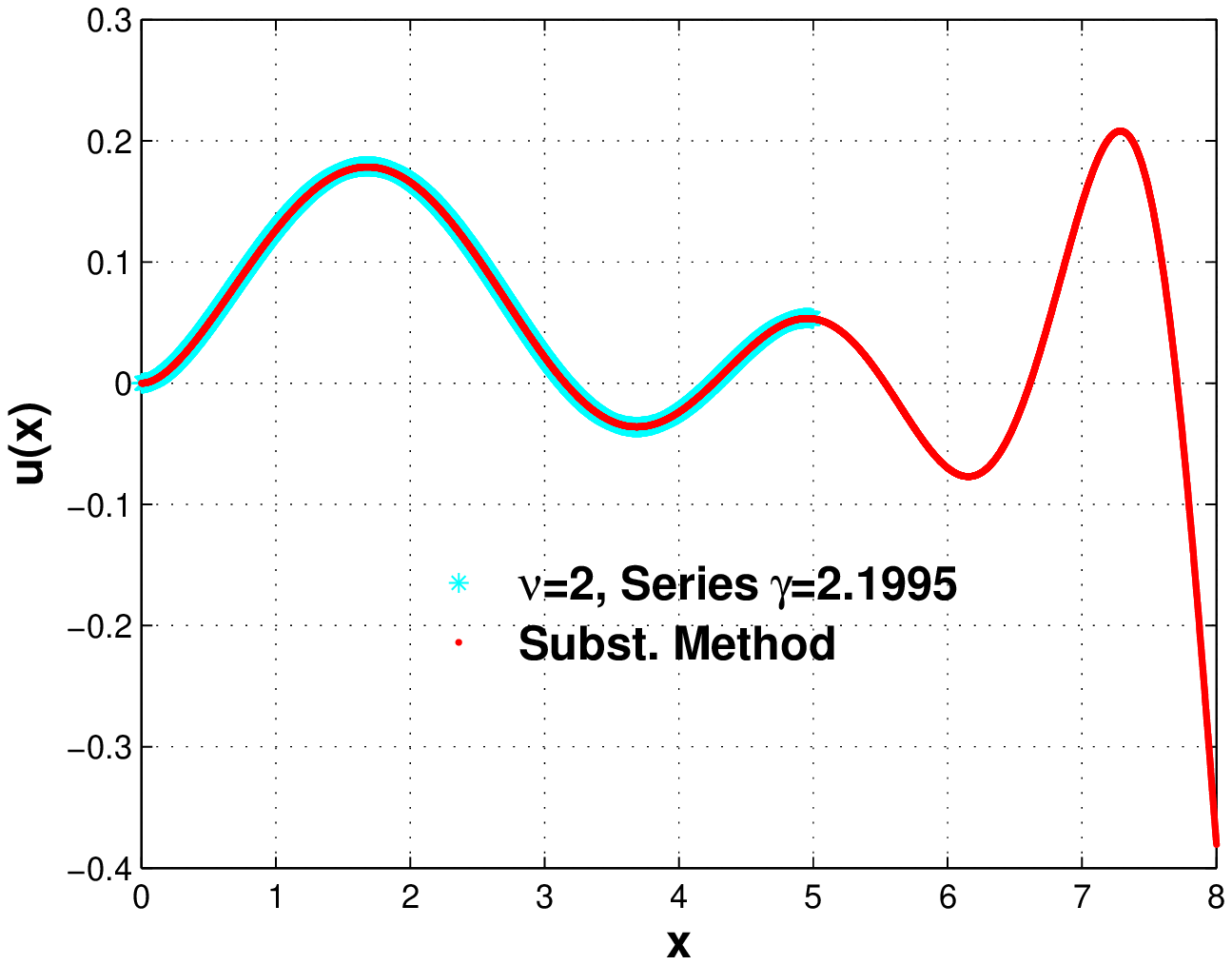}
		\includegraphics[width=6cm]{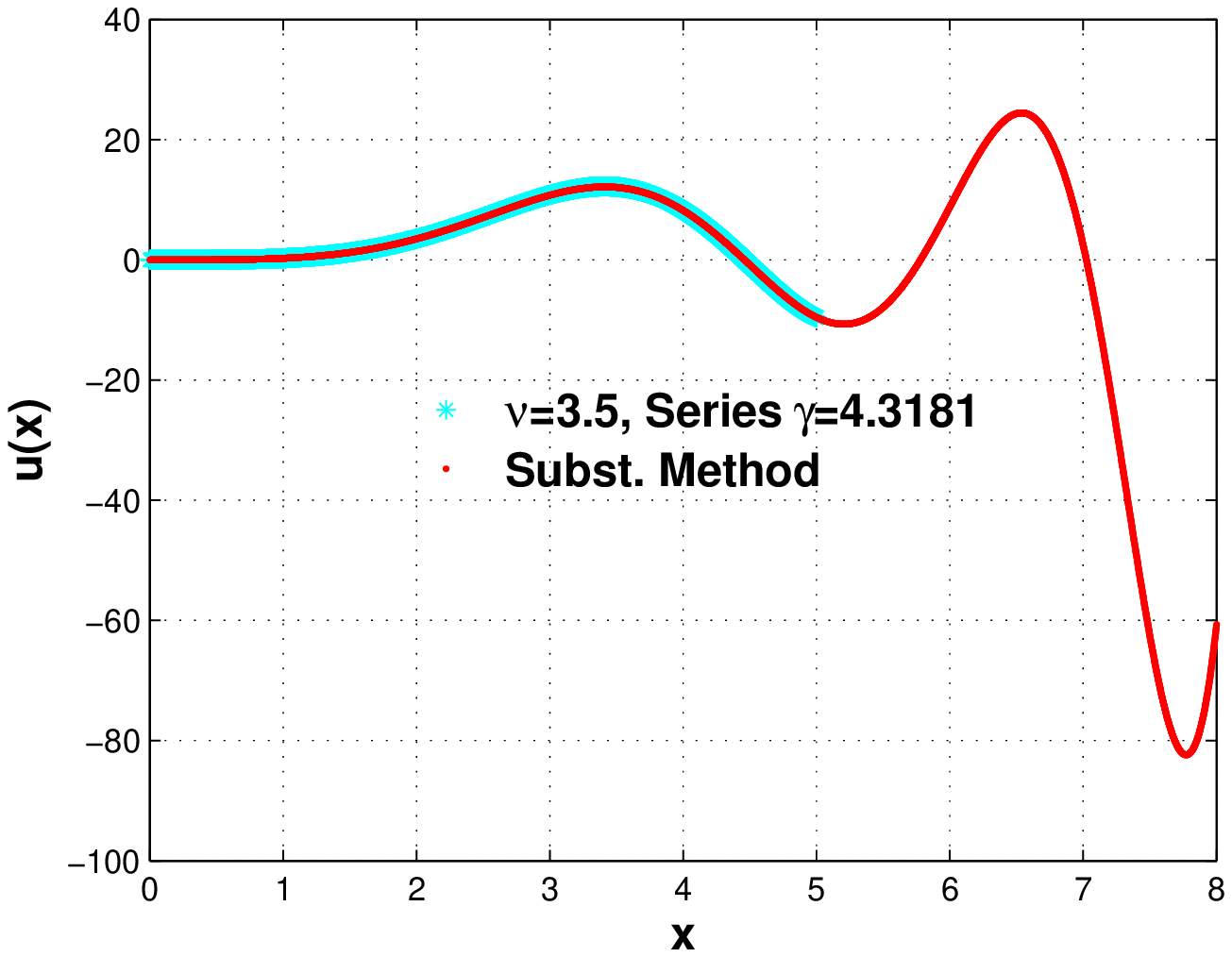}
		\caption{Example 1. Light Blue - exact series solution, red - solution generated by the substitution method \cite{DS2020DuSl-2}. On [0,5] graphs of both methods match.}\label{SolEx1Subst}	
	\end{center}
\end{figure}

The solution becomes unstable on a larger interval because $x^{\gamma+s n}$ becomes a very large number and more terms in the sequence are needed for accurate results but these terms are a product of a very small coefficient $c_n$ and very large power of $x$ which numerically becomes unstable.  If a solution is needed on a larger interval we can solve the equation using the substitution technique \cite{DS2020DuSl-2} after making sure that it matches on a smaller interval. The solution at zero is zero and its first derivative is also zero.  With these initial conditions the substitution method only finds the trivial solution.  We can slightly perturb the system by assigning a small value for the derivative and, as we know, it only effects the constant in front of the function $u(x)$.  After this minor perturbation/calibration of the first derivative we can generate the solution on a larger interval.  These solutions are shown in Figure \ref{SolEx1Subst}.
\qed \\
\\
{\scshape Example 2 } (constant coefficients, integer derivatives).\\

 This example demonstrates how the methodology works on a very simple example for which we know the analytic solution.  
 \begin{equation*}
 u'+u=0.
 \end{equation*}
 We know that its solution is $u=Ce^{-x}$.  Our method converts this equation into
 \begin{equation*}
 xu'+xu=0
 \end{equation*}
 and states that series solution exists and is in the form $u(x)=\ds\sum_{n=0}^{\infty}c_n x^{sn}$.  \\
 Solution for characteristic equation 
 \begin{equation*}
 \frac{\Gamma(1+\gamma)}{\Gamma(1+\gamma-1)}=0
 \end{equation*}
 is $\gamma=0$.  Based on \eqref{BetasNs} step $s=1$. Therefore, we get
 \begin{equation*}
 c_n=c_n^\beta=-\frac{c_{n-1}}{Q(n,1)}=-c_{n-1}\frac{\Gamma(1+n-1)}{\Gamma(1+n)}=-\frac{c_{n-1}}{n}=c_0\frac{(-1)^n}{n!}
 \end{equation*}
 and the solution as expected is
 \begin{equation*}
 u(x)=\sum_{n=0}^{\infty}c_n x^{sn}=c_0\sum_{n=0}^{\infty}\frac{(-1)^n}{n!}x^n=c_0 e^{-x}.
 \end{equation*} \qed
\\
   \\
{\scshape Example 3 } (fractional derivative with constant coefficients).\\

Let us consider equation 4.1.21 \cite{Kilbas} with specific values of parameters 
\begin{eqnarray}\label{Example1Kilb}
D_R^{1.7}u(x)-2u(x)=0 
\end{eqnarray}
with initial conditions
\begin{eqnarray}\label{Ex1IC}
D_R^{0.7}u(0)=1.2, D_R^{-0.3}u(0)=1.5. 
\end{eqnarray}
%
%
The characteristic equation for \eqref{Example1Kilb} becomes
\begin{equation}
\frac{\Gamma(1+\gamma)}{\Gamma(1+\gamma-1.7)}=0.
\end{equation}
It has two roots with $\gamma > -1: \gamma_1=0.7, \gamma_2=-0.3$.  Step $s=1.7$ \\
 (no $p_i; \beta=1.7$).  Therefore we get two series solutions:
\begin{eqnarray}
u(x)&=&\sum_{n=0}^{\infty}c_n^1 x^{0.7+1.7n}, \nonumber \\
u(x)&=&\sum_{n=0}^{\infty}c_n^2 x^{-0.3+1.7n}. \nonumber
\end{eqnarray}
In this case (since coefficient in front of $u(x)$ needs to be one) $d_1=-0.5$ and therefore
\begin{eqnarray}
c_n^1=-\frac{c_{n-1}^1\Gamma(1+0.7+1.7n-1.7)}{-0.5\Gamma(1+0.7+1.7n)}=c_{n-1}^1\frac{2\Gamma(1.7n)}{\Gamma(1.7+1.7n)}=c_0^1\frac{2^n\Gamma(1.7)}{\Gamma(1.7+1.7n)} \ \nonumber \\
c_n^2=-\frac{c_{n-1}^2\Gamma(1-0.3+1.7n-1.7)}{-0.5\Gamma(1-0.3+1.7n)}=c_{n-1}^2\frac{2\Gamma(1.7n-1)}{\Gamma(0.7+1.7n)}=c_0^2\frac{2^n\Gamma(0.7)}{\Gamma(0.7+1.7n)}.\nonumber 
\end{eqnarray}
Then, the solution is
\begin{eqnarray}\label{Ex1Soln}
u(x)&=&c_0^1\sum_{n=0}^{\infty}\frac{\Gamma(1.7) 2^n x^{ 0.7+1.7n}}{\Gamma(1.7+1.7n)} + 
c_0^2\sum_{n=0}^{\infty}\frac{\Gamma(0.7) 2^n x^{-0.3+1.7n}}{\Gamma(0.7+1.7n)} \nonumber \\ 
&=&c_0^1\Gamma(1.7) \sum_{n=0}^{\infty}\frac{2^n x^{ 0.7+1.7n}}{\Gamma(1.7+1.7n)} + 
c_0^2\Gamma(0.7) \sum_{n=0}^{\infty}\frac{2^n x^{-0.3+1.7n}}{\Gamma(0.7+1.7n)}. \ \ \ \: \:  
\end{eqnarray}
Let us see what initial condition \eqref{Ex1IC} gives us \cite{Podlubny}:
\begin{eqnarray}
D_R^{0.7}u(0)&=&c_0^1\Gamma(1.7)\frac{\Gamma(1+0.7)}{\Gamma(1.7)\Gamma(1+0.7-0.7)}+0=c_0^1 \Gamma(1.7)=1.2 \nonumber \\
D_R^{-0.3}u(0)&=&0+c_0^2\Gamma(0.7)\frac{\Gamma(1-0.3)}{\Gamma(1-0.3+0.3)\Gamma(0.7)}=c_0^2 \Gamma(0.7)=1.5, \nonumber
\end{eqnarray} 
which makes solution \eqref{Ex1Soln} the same as previously identified in \cite{Kilbas} analytical solution. \qed	\\
\\
{\scshape Example 4 } (fractional derivatives with powers of $x$).\\

Let us apply the methodology to another example from  \cite{Kilbas}, Example 4.1.36.
\begin{eqnarray}\label{ExXpower}
D_R^{0.5}u(x)-\lambda x^\beta u(x)=0, D_R^{-0.5}u(0)=b.
\end{eqnarray}
For simplicity sake we will assume that $\beta=0.7 \in \mathbb{Q^+}$.  Then, our interpretation of this problem is
\begin{eqnarray}\label{ExXpower2}
-\frac{1}{\lambda}x^{0.5}D_R^{0.5}u(x)+x^{1.2} u(x)=0, D_R^{-0.5}u(0)=b.
\end{eqnarray}
Characteristic equation \eqref{CondC0ConstCoeffs} produces one root $\gamma=-0.5$, step $s=1.2$, and $n_\beta=1$.
Then, the solution to equation \eqref{ExXpower2} is  
\begin{equation}\label{ExXposerSoln}
u(x)=\sum_{n=0}^{\infty}c_n x^{-0.5+1.2n}, 
\end{equation}
where coefficients are expressed as
\begin{eqnarray}
c_n&=&-\frac{-\lambda c_{n-1}^1\Gamma(1-0.5+1.2n-0.5)}{\Gamma(1-0.5+1.2n)} =c_{n-1}\frac{\lambda\Gamma(1.2n)}{\Gamma(-0.5+1.2n)} \nonumber \\
&=&c_0\lambda^n\prod_{j=0}^{n-1}\frac{\Gamma(1.2+1.2n)}{\Gamma(1.7+1.2n)} \nonumber 
\end{eqnarray}
and the solution then is 
\begin{eqnarray}\label{ExXposerSolnF}
u(x)&=&c_0\sum_{n=0}^{\infty}\lambda^n\prod_{j=0}^{n-1}\frac{\Gamma(1.2+1.2n)}{\Gamma(1.7+1.2n)} x^{-0.5+1.2n}\nonumber \\
&=&c_0x^{-0.5}\sum_{n=0}^{\infty}\prod_{j=0}^{n-1}\frac{\Gamma(1.2+1.2n)}{\Gamma(1.7+1.2n)} \lambda^n x^{1.2n}
\end{eqnarray}
If we apply the initial condition then we obtain
\begin{equation*}
D_R^{-0.5}u(0)=c_0\frac{\Gamma(0.5)}{\Gamma(1)}=c_0\Gamma(0.5)=c_0\sqrt{\pi}.
\end{equation*}
Then, $c_0=\ds\frac{b}{\sqrt{\pi}}$ and we arrive at solution 
\begin{equation}\label{ExXposerSolnF2}
u(x)=\frac{b}{\sqrt{\pi}}x^{-0.5}\sum_{n=0}^{\infty}\prod_{j=0}^{n-1}\frac{\Gamma(1.2+1.2n)}{\Gamma(1.7+1.2n)} \lambda^n x^{1.2n}=\frac{b}{\sqrt{\pi}}x^{-0.5}E_{0.5,2.4,0.4}[\lambda x^{1.2}],
\end{equation}
where $E_{\alpha,m,l}(z)=\ds\sum_{k=0}^{\infty}c_k z^k$ with $c_0=1, c_k=\ds\prod_{j=0}^{k-1}\frac{\Gamma[\alpha(jm+l)+1]}{\Gamma[\alpha(jm+l+1)+1]}$. $E_{\alpha,m,l}(z)$ is the generalized Mittag-Leffler function introduced by Kilbas and Saigo.  Solution \eqref{ExXposerSolnF2} matches the solution identified in \cite{Kilbas} for this problem. \qed \\

\section{Conclusions}
The quasi-Bessel equations have terms like $x^{\alpha + p}D^\alpha u(x)$ and extend further the multi-term 
 fractional Bessel equations and cover an essentially broader class of equations. Now the powers of $x$ may avoid matching the order of corresponding fractional derivatives (except for the highest derivative), whereas such a match is usually required for the successful applications of integral Mellin transform. Particularly, the Cauchy-Euler and constant-coefficient equations are included in the class of quasi-Bessel equations.
 
If the deviating (shifting) parameters $p$ are non-negative and meet certain non-restrictive conditions, then we construct the existence theory for quasi-Bessel equations in the class of fractional series solutions. Also, we prove the uniqueness theorem for the initial value problem.

Several presented examples and computational experiments, including the examples for constant-coefficient equations, reaffirm the validity of our methodology.

The quasi-Bessel equations with $p < 0$ do not fit the presented theory and
remain an open problem.

\section*{Acknowledgments}
The authors would like to thank Professor L.~Boyadjiev for drawing our attention to the fractional Bessel equation.
Also, the authors are obliged to Professor Kiryakova for valuable recommendations. 

\bigskip \smallskip

\it

\noindent
$^1$ Stevens Institute of Technology,
1 Castle Point Terrace\\
Hoboken, NJ 07030, USA  \\[4pt]
e-mail: pdubovsk@stevens.edu, jslepoi@stevens.edu
\hfill  \\[12pt]

\
\end{document}